\documentclass{amsart}

\usepackage{array,amsmath,amsthm,amsxtra,amssymb,graphicx,subfigure}

\usepackage{supertabular}
\usepackage[absolute]{textpos}

\textblockorigin{0in}{0in}
\newtheorem{theorem}{Theorem}
\newtheorem{corollary}[theorem]{Corollary}
\newtheorem{lemma}[theorem]{Lemma}
\newtheorem{prop}[theorem]{Proposition}

\theoremstyle{definition}
\newtheorem{defn}[theorem]{Definition}

\newtheorem{remark}[theorem]{Remark}

\numberwithin{theorem}{section}

\newcommand{\hg}[4]{
_{2}F_{1} \left( 
\begin{matrix}
#1, & #2 \\
   & #3 \\
\end{matrix}
\bigg{\vert} #4
\right) 
}

\newcommand{\fp}
{\mathbb{F}_p}

\newcommand{\fpc}
{\mathbb{F}_p^{\times}}

\newcommand{\fphat}
{\widehat{\mathbb{F}_{p}^{\times}}}

\newcommand{\C}{\mathbb{C}}

\newcommand{\Z}{\mathbb{Z}}

\begin{document}

\pagenumbering{arabic}
\setlength{\headheight}{12pt}
\pagestyle{myheadings}

\title[Hypergeometric Functions over $\fp$]{Hypergeometric functions over $\fp$ and relations to elliptic curves and modular forms}

\author{Jenny G. Fuselier}
\maketitle

\date{}
\maketitle

\begin{abstract}
For primes $p\equiv 1 \pmod{12}$, we present an explicit relation between the traces of Frobenius on a family of elliptic curves with $j$-invariant $\frac{1728}{t}$ and values of a particular $_2F_1$-hypergeometric function over $\fp$.  Additionally, we determine a formula for traces of Hecke operators $\textnormal{T}_k(p)$ on spaces of cusp forms of weight $k$ and level 1 in terms of the same traces of Frobenius.  This leads to formulas for Ramanujan's $\tau$-function in terms of hypergeometric functions.
\end{abstract}

%%%%%%%%%%%%%%%%%%%%%%%%%%%%%%%%%%%%%%%%
%% Section 1 - Intro and Main Results %%
%%%%%%%%%%%%%%%%%%%%%%%%%%%%%%%%%%%%%%%%

\section{Introduction and Statement of Main Results}\label{intro}

Let $p$ be a prime and let $\widehat{\mathbb{F}_{p}^{\times}}$ denote the group of all multiplicative characters on $\fpc$.   We extend $\chi\in\ \widehat{\mathbb{F}_{p}^{\times}}$ to all of $\fp$ by setting $\chi(0)=0$.   For $A,B\in\fphat$, let $J(A,B)$ denote the usual Jacobi symbol and define  

\begin{equation}\label{binom coeff}
\binom{A}{B}:=\frac{B(-1)}{p}J(A,\overline{B})=\frac{B(-1)}{p} \sum_{x\in \fp} A(x)\overline{B}(1-x).
\end{equation}

\noindent In the 1980s, Greene \cite{Gr87} defined \emph{hypergeometric functions over $\fp$} in the following way:   

\begin{defn}[\cite{Gr87} Defn. 3.10]\label{hg}
If $n$ is a positive integer, $x\in\fp$, and $A_0,A_1,\dots,A_n,$\\$B_1,B_2,\dots,B_n \in \widehat{\mathbb{F}_{p}^{\times}}$, then define
$$ _{n+1}F_{n} \left( 
\begin{matrix}
A_0, & A_1, & \dots, & A_n \\
     & B_1, & \dots, & B_n \\
\end{matrix}
\bigg{\vert} x
\right) 
:= \frac{p}{p-1} \sum_{\chi\in\widehat{\mathbb{F}_{p}^{\times}}} \binom{A_0\chi}{\chi} \binom{A_1\chi}{B_1\chi} \dots \binom{A_n\chi}{B_n\chi} \chi(x).$$
\end{defn}

Greene explored the properties of these functions and showed that they satisfy many transformations analogous to those enjoyed by their classical counterparts.  The introduction of these hypergeometric functions over $\fp$ generated interest in finding connections they may have with modular forms and elliptic curves.  In recent years, many results have been proved in this direction, and the main results given here are similar in nature. 

Throughout, we consider a family of elliptic curves having $j$-invariant $\frac{1728}{t}$.  Specifically, for $t\in\fp$, $t\neq0,1 $ we let 
\begin{equation}\label{defn of E}
E_t: y^2=4x^3-\frac{27}{1-t}x - \frac{27}{1-t}.
\end{equation}
Further, we let $a(t,p)$ denote the trace of the Frobenius endomorphism on $E_t$.  In particular, for $t\neq0,1$, we have
$$a(t,p)=p+1-\#E_t(\fp),$$ where $\#E_t(\fp)$ counts the number of solutions to $y^2 \equiv 4x^3-\frac{27}{1-t}x - \frac{27}{1-t}$ (mod $p$), including the point at infinity.

Henceforth, we let $p$ be a prime number with $p\equiv 1$ (mod 12).  With this in mind, we let $\xi\in \widehat{\mathbb{F}_{p}^{\times}}$ have order 12.  Also, we denote by $\varepsilon$ and $\phi$ the trivial and quadratic characters, respectively.  In this setting, our first main result explicitly relates the above trace of Frobenius and the values of a hypergeometric function over $\fp$.

\begin{theorem}\label{first main theorem}
Suppose $p\equiv 1 \pmod{12}$ is prime and $\xi\in \widehat{\mathbb{F}_{p}^{\times}}$ has order $12$.  Then, if $t\in\fp\backslash\{0,1\}$ and notation is as above, we have
$$p\,\hg{\xi}{\xi^5}{\varepsilon}{t}= \psi(t) a(t,p),$$
where $\psi(t)= -\phi(2)\xi^{-3}(1-t)$.
\end{theorem}

After setting up the necessary preliminaries, we give the proof of Theorem \ref{first main theorem} in Section \ref{proof of theorem 1}.  The second main result utilizes the same family of elliptic curves, $E_t$, to obtain a trace formula for Hecke operators on spaces of cusp forms in level $1$.

Let $\Gamma=SL_2(\Z)$ and let $M_k$ and $S_k$, respectively, denote the spaces of modular forms and cusp forms of weight $k$ for $\Gamma.$  Further, let $\textnormal{Tr}_k(\Gamma,p)$ denote the trace of the Hecke operator $\textnormal{T}_k(p)$ on $S_k$.  Then the following completely describes these traces $\textnormal{Tr}_k(\Gamma,p)$, for a certain class of primes $p$:

\begin{theorem}\label{second main theorem}
Suppose $p\equiv 1 \pmod{12}$ is prime.  Let $a,b\in\mathbb{Z}$ such that  $p=a^2+b^2$ and $a+bi \equiv 1 \,(2+2i)$ in $\mathbb{Z}[i]$.  Also, let $c,d\in\mathbb{Z}$ such that $p=c^2-cd+d^2$ and $c+d\omega \equiv 2 \,(3)$ in $\mathbb{Z}[\omega]$, where $\omega=e^{2\pi i/3}$.  Then for even $k\geq4$,   
$$ \textnormal{Tr}_k(\Gamma,p)=-1-\lambda(k,p)-\sum_{t=2}^{p-1} G_k(a(t,p),p),$$ where
$$\lambda(k,p)=\frac{1}{2}[G_k(2a,p)+G_k(2b,p)]+\frac{1}{3}[G_k(c+d,p)+G_k(2c-d,p)+G_k(c-2d,p)]$$ and
$$G_k(s,p)=\sum_{j=0}^{\frac{k}{2}-1} (-1)^j \binom{k-2-j}{j} p^j s^{k-2j-2}.$$
\end{theorem}

Combining Theorems \ref{first main theorem} and \ref{second main theorem} gives a way of writing the traces $\textnormal{Tr}_k(\Gamma,p)$ in terms of hypergeometric functions.  Formulas for Ramanujan's $\tau$-function follow by taking $k=12$ in Theorem \ref{second main theorem}.

Finally, Theorem \ref{second main theorem} gives rise to an inductive formula for the traces $\textnormal{Tr}_k(\Gamma,p)$, in terms of hypergeometric functions.  To state it, we utilize the notation for $G_k(s,p)$ and $\lambda(k,p)$ given in Theorem \ref{second main theorem}.  

\begin{theorem}\label{recursion}
Suppose $p\equiv 1 \pmod{12}$ is prime. Let $k\geq 4$ be even and define $m=\frac{k}{2}-1$.  Then
\begin{align*}
\textnormal{Tr}_{2(m+1)}(\Gamma,p)&=-1-\lambda(2m+2,p)+b_0(p-2)-\sum_{t=2}^{p-1}p^{2m}\phi^m(1-t)\hg{\xi}{\xi^5}{\varepsilon}{t}^{2m}\\
&\quad\,-\sum_{i=1}^{m-1}b_i(1+\lambda(2i+2,p))-\sum_{i=1}^{m-1}b_i\textnormal{Tr}_{2i+2}(\Gamma,p),
\end{align*}
where
$$b_i=p^{m-i}\left[\binom{2m}{m-i}-\binom{2m}{m-i-1}\right].$$

\end{theorem}

In Section \ref{history hg and ec}, we recall a few recent results that relate counting points on varieties over $\fp$ to hypergeometric functions and make comparisons with classical hypergeometric functions.  In Section \ref{preliminaries}, we introduce necessary preliminaries and give the proof of Theorem \ref{first main theorem} in Section \ref{proof of theorem 1}.  In Sections \ref{history hg and mf} and \ref{proof of theorem 2}, we focus on Theorem \ref{second main theorem}, beginning with a history of similar results, and building to a proof of the theorem.  This is followed by Section \ref{recursion proof}, in which we give a proof of Theorem \ref{recursion}.  We close in Section \ref{tau corollaries} with various corollaries relating Ramanujan's $\tau$-function to hypergeometric functions.  Specifically, Corollary \ref{tau corollary 2} expresses $\tau(p)$ explicitly in terms of tenth powers of a $_2F_1$-hypergeometric function.

%%% STILL NEED TO REWORK THE ABOVE PARAGRAPH AFTER I KNOW WHAT THE FINAL SECTIONING WILL BE FOR THE MF STUFF%%%%%

\section{Recent History: Hypergeometric Functions and Elliptic Curves}\label{history hg and ec}
%Need to include history about connections between hgf's and elliptic curves. Need the following:

%\begin{itemize}

%\item stiller's result, explicitly stated (POSTPONE?)
%\item beukers results (maybe not stated explicity??) (DONE)
%\item koike and ono's results (DONE)
%\item maybe matt and sharon's result (published??) connecting counting points on other varieties with values of certain hgfs

%\end{itemize}

Relationships between hypergeometric functions over $\fp$ and elliptic curves are perhaps not surprising, as classical hypergeometric series have many known connections to elliptic curves.  An important example of these classical series is defined for $a,b,c \in \mathbb{C}$ as 
$$_2F_1[a,b;c;z]:=\sum_{n=1}^{\infty} \frac{(a)_n(b)_n}{(c)_nn!}z^n,$$
where $(w)_n=w(w+1)(w+2)\cdots (w+n-1)$.

The specialization $_2F_1[\frac{1}{2},\frac{1}{2};1;t]$ is related to the Legendre family of elliptic curves, as it is a constant multiple of an elliptic integral which represents a period of the associated lattice.  More recently, Beukers \cite{Be93} gave identifications between periods of families of elliptic curves and values of particular hypergeometric series.  For example, he related a period of $y^2=x^3+tx+1$ to the values $_2F_1[\frac{1}{12},\frac{7}{12};\frac{2}{3};-\frac{4}{27}t^3]$ and a period of $y^2=x^3-x-t$ to the values $_2F_1[\frac{1}{12},\frac{5}{12};\frac{1}{2};\frac{27}{4}t^2]$.  We give finite field analogues of these results at the close of the next section.

Following Greene's introduction of hypergeometric functions over $\fp$ in the 1980s, results emerged linking their values to counting points on varieties over $\fp$.  Let $\phi$ and $\varepsilon$ denote the unique quadratic and trivial characters, respectively, on $\fpc$.  Further, define two families of elliptic curves over $\fp$ by
\begin{align*}
_2E_1(t):y^2=x(x-1)(x-t)\\
_3E_2(t):y^2=(x-1)(x^2+t).
\end{align*}
Then, for odd primes $p$ and $t\in\fp$, define the traces of Frobenius on the above families by
\begin{align*}
_2A_1(p,t)&=p+1- \#_2E_1(t)(\fp),\quad t\neq0,1\\
_3A_2(p,t)&=p+1- \#_3E_2(t)(\fp),\quad t\neq0,-1.
\end{align*}
These families of elliptic curves are closely related to particular hypergeometric functions over $\fp$.  For example, $\hg{\phi}{\phi}{\varepsilon}{t}$ arises in the formula for Fourier coefficients of a modular form associated to $_2E_1(t)$ (\cite{Ko92,On98}).  Further, Koike and Ono, respectively, gave the following explicit relationships:

\begin{theorem}[(a) Koike \cite{Ko92}, (b) Ono \cite{On98}]\label{koike ono}
Let $p$ be an odd prime.  Then

$(a)\quad p \, \hg{\phi}{\phi}{\varepsilon}{t}=-\phi(-1) _2A_1(p,t), \quad t\neq0,1$

$(b)\quad p^2 \, _{3}F_{2} \left( 
\begin{matrix}
\phi, & \phi, & \phi \\
     & \varepsilon, & \varepsilon \\
\end{matrix}
\bigg{\vert} 1+\frac{1}{t}
\right) 
=
\phi(-t)(_3A_2(p,t)^2-p), \quad t\neq 0,-1.$
\end{theorem}

%State Sharon and Matt's result here too (find out if it's published, so I know if I can cite it here).

%%%%%%%%%%%%%%%%%%%%%%%%%%%
% Preliminaries for Thm 1 %
%%%%%%%%%%%%%%%%%%%%%%%%%%%

\section{Preliminaries on Characters and Hypergeometric Functions}\label{preliminaries}

The proof of Theorem \ref{first main theorem} involves two main steps.  First, we derive a formula for $a(t,p)$ in terms of Gauss sums, and then we write the hypergeometric function in terms of Gauss sums.  The final proof follows from comparing the two.  Before doing that, we fix notation and recall some basic facts regarding Gauss sums.

Throughout, let $p\equiv 1 \pmod{12}$ be prime.  If $A\in\fphat$, let 
$$ G(A)=\sum_{x\in\fp}A(x)\zeta^{x}$$
denote the Gauss sum, with $\zeta=e^{2\pi i/p}$.

Since $\fpc$ is cyclic, let $T$ denote a fixed generator of the character group, i.e. $\langle T \rangle= \widehat{\mathbb{F}_{p}^{\times}}$.  With this in mind, we often use the notation $G_m:=G(T^m)$.  Recall the following elementary properties of Gauss sums, the proofs of which can be found in Chapter 8 of \cite{IR90}:

\begin{lemma}\label{orthog relation}
Let $T$ be a generator for $\widehat{\mathbb{F}_{p}^{\times}}$.  Then

$(a) \quad \displaystyle{\sum_{x\in\fp}T^n(x)=
\begin{cases}
p-1 & \textnormal{if} \,\, T^n=\varepsilon\\
0 & \textnormal{if} \,\, T^n\neq\varepsilon
\end{cases}}$

$(b) \quad \displaystyle{\sum_{n=0}^{p-2}T^n(x)=
\begin{cases}
p-1 &  \textnormal{if} \,\, x=1\\
0 & \textnormal{if} \,\, x \neq 1.
\end{cases}}$

\end{lemma}

The next result calculates the values of two particular Gauss sums, $G(\varepsilon)$ and $G(\phi)$.

\begin{lemma}\label{special Gauss sums}
$(a) \quad G(\varepsilon)=G_0=-1$ \\
\hspace*{0.9in}$(b) \quad \displaystyle{G(\phi)=G_\frac{p-1}{2}=
\begin{cases} \sqrt{p} & \textnormal{if}\,\, p\equiv 1 \pmod{4}\\
i\sqrt{p} & \textnormal{if}\,\,p\equiv 3 \pmod{4}.
\end{cases}}$
\end{lemma}

We now define an additive character 

\begin{eqnarray*}\label{theta def}
\theta:\fp  &\rightarrow\C \\ 
\theta(\alpha)&=\zeta^{\alpha}. 
\end{eqnarray*}

Notice that we can write Gauss sums in terms of $\theta$, as we have $G(A)=\sum_{x\in\fp} A(x)\theta(x)$.  In addition, the following lemma describing $\theta$ in terms of Gauss sums is straightforward to prove via the orthogonality relations given above:

\begin{lemma}\label{theta equality} For all $\alpha\in\fpc$,
$$\theta(\alpha)=\frac{1}{p-1} \sum_{m=0}^{p-2} G_{-m}T^m(\alpha).$$
\end{lemma}

We also require a few properties of hypergeometric functions over $\fp$ that Greene proved in \cite{Gr87}.  The first provides a formula for the multiplicative inverse of a Gauss sum.

\begin{lemma}[\cite{Gr87} Eqn. 1.12]\label{greene 112}
If $k\in\mathbb{Z}$ and $T^k\neq\varepsilon$, then 
$$G_kG_{-k}=pT^k(-1).$$
\end{lemma}

The following result was given by Greene as the definition of the hypergeometric function when $n=1$.  It provides an alternative to Definition \ref{hg}, and in particular, it allows us to write the $_2F_1$ hypergeometric function as a character sum.
\begin{theorem}[\cite{Gr87} Defn. 3.5]\label{greene 35}
If $A,B,C \in \widehat{\mathbb{F}_{p}^{\times}}$ and $x \in \fp$, then $$\hg{A}{B}{C}{x} = \varepsilon(x) \frac{BC(-1)}{p} \sum_{y=0}^{p-1} B(y)\overline{B}C(1-y)\overline{A}(1-xy).$$
\end{theorem}

In \cite{Gr87}, Greene presented many transformation identities satisfied by the hypergeometric functions he defined.  The theorem below allows for the argument $x\in\fp$ to be replaced by $1-x$.

\begin{theorem}[\cite{Gr87} Theorem 4.4]\label{greene 44}
If  $A,B,C \in \widehat{\mathbb{F}_{p}^{\times}}$ and $x\in\fp \backslash \{0,1\}$, then
$$\hg{A}{B}{C}{x}=A(-1)\hg{A}{B}{AB\overline{C}}{1-x}.$$
\end{theorem}

Next,we recall a classical relationship between Gauss and Jacobi sums, but we write it utilizing Greene's definition for the binomial coefficient, given in \eqref{binom coeff}.

\begin{lemma}\label{greene 29}
If $T^{m-n}\neq \varepsilon$, then
$$ \binom{T^m}{T^n}=\frac{G_mG_{-n}T^n(-1)}{G_{m-n} \cdot p}.$$ 
\end{lemma}
% Maybe move this "`easy"' lemma earlier?

%%%%%%%%%%%%%%%%%%%%%%%%%%%%%%%%%%%%%%%%%%%
% Hasse Davenport Relation and corollaries%
%%%%%%%%%%%%%%%%%%%%%%%%%%%%%%%%%%%%%%%%%%%

The final relation on characters that is necessary for our proof is the \emph{Hasse-Davenport relation.}  The most general version of this relation involves an arbitrary additive character, and can be found in \cite{La90}.  We require only the case when $\theta$ is taken as the additive character:

\begin{theorem}[Hasse-Davenport Relation \cite{La90}]\label{Hasse Davenport}
Let $m$ be a positive integer and let $p$ be a prime so that $p\equiv1 \pmod{m}.$  Let $\theta$ be the additive character on $\mathbb{F}_p$ defined by $\theta(\alpha)=\zeta^{\alpha}$, where $\zeta=e^{2\pi i/p}$.  For multiplicative characters $\chi,\psi\in \widehat{\mathbb{F}_{p}^{\times}}$, we have
$$\prod_{\chi^m=1}G(\chi\psi)=-G(\psi^m)\psi(m^{-m})\prod_{\chi^m=1}G(\chi).$$
\end{theorem}

\begin{proof}
See \cite{La90}, page 61.
\end{proof}

The proof of Theorem \ref{first main theorem} requires two special cases of the Hasse-Davenport relation, which are easily verified by taking $m=2$ and $m=3$, respectively.

\begin{corollary}\label{HD2}  If $p\equiv 1 \pmod{4}$ and $k\in\mathbb{Z}$,
$$G_{-k}G_{-\frac{p-1}{2}-k}=\sqrt{p}\,G_{-2k}T^k(4).$$
\end{corollary}

\begin{corollary}\label{HD3} If $k\in\mathbb{Z}$ and $p$ is a prime with $p\equiv 1 \pmod{3}$ then
$$G_kG_{k+\frac{p-1}{3}}G_{k+\frac{2(p-1)}{3}}=p\,T^{-k}(27)T^{\frac{p-1}{3}}(-1)G_{3k}.$$
\end{corollary}

%%%%%%%%%%%%%%%%%%%%%%%%%%%%%%%
% Proof of First Main Theorem %
%%%%%%%%%%%%%%%%%%%%%%%%%%%%%%%

\section{Proof of Theorem \ref{first main theorem}}\label{proof of theorem 1}

We begin by deriving a formula for the trace of Frobenius in terms of Gauss sums.  Throughout this section, let $s=\frac{p-1}{12}$ and define $P(x,y)= y^2-4x^3+\frac{27}{1-t}x+ \frac{27}{1-t}$.

%\begin{notation}
%$\bullet$ Let $s=\frac{p-1}{12}.$\\
%\hspace*{0.735in}$\bullet$ Let  $P(x,y)= y^2-4x^3+\frac{27}{1-t}x+ \frac{27}{1-t}$.
%\end{notation}

Recall from the previous section that $\theta$ is the additive character on $\fp$ given by $\theta(\alpha)=\zeta^{\alpha}$, where $\zeta=e^{2\pi i/p}$.  Note that if $(x,y)\in\fp^2$, then $$\sum_{z\in\fp}\theta(zP(x,y))=\begin{cases} p & \textnormal{if} \,P(x,y)=0\\ 0 &\textnormal{if}\, P(x,y)\neq 0. \end{cases}$$  So we have 
\begin{align*}
p\cdot(\#E_t(\fp)-1) &= \sum_{z\in\fp} \sum_{x,y \in\fp} \theta(zP(x,y))\\ 
        &= \sum_{x,y,\in\fp} 1 + \sum_{z\in\fpc}\sum_{x,y\in\fp}\theta(zP(x,y)),
\end{align*}
after breaking apart the $z=0$ contribution.  Then, by separating the sums according to whether $x$ and $y$ are $0$ and applying the additivity of $\theta$, we have
\begin{align*}
p\cdot(\#E_t(\fp)-1)&= p^2+ \sum_{z\in\fpc}\theta\left(z\frac{27}{1-t}\right)+\sum_{z\in\fpc}\sum_{y\in\fpc} \theta(zy^2)\theta\left(z\frac{27}{1-t}\right)\\
        &\hspace*{0.2in} + \sum_{z\in\fpc}\sum_{x\in\fpc}\theta(-4zx^3)\theta\left(zx\frac{27}{1-t}\right)\theta\left(z\frac{27}{1-t}\right)\\
        &\hspace*{0.2in}+\sum_{x,y,z\in\fpc}\theta(zP(x,y))\\
        &:= p^2+A+B+C+D,
\end{align*}
where $A$, $B$, $C$, and $D$ are set to be the four sums appearing in the previous line. These four sums are computed using Lemmas \ref{orthog relation}, \ref{special Gauss sums} and \ref{theta equality} repeatedly. We provide the computation for $D$ here, which requires the most steps.  The other three follow in a similar manner.  We begin with four applications of Lemma \ref{theta equality} and find that

%%%%%%%%%%%%%%%%%%%
%% D Computation %%
%%%%%%%%%%%%%%%%%%%

\allowdisplaybreaks{
\begin{align*}
D &=  \frac{1}{(p-1)^4}\sum_{x,y,z \in\fpc}\sum_{j,k,\ell,m=0}^{p-2}G_{-j}G_{-k}G_{-\ell}G_{-m}T^j(zy^2)T^k(-4zx^3)\\&\qquad\qquad\qquad \cdot T^{\ell}\left(zx\frac{27}{1-t}\right)T^m\left(z\frac{27}{1-t}\right)\\ 
  &=  \frac{1}{(p-1)^4}\sum_{x,y \in\fpc}\sum_{j,k,\ell,m=0}^{p-2}G_{-j}G_{-k}G_{-\ell}G_{-m}T^j(y^2)T^k(-4x^3)T^{\ell}\left(x\frac{27}{1-t}\right)\\
  & \qquad \qquad \qquad \cdot T^m\left(\frac{27}{1-t}\right)\sum_{z\in\fpc}T^{j+k+\ell+m}(z),
\end{align*}
} 
\noindent after simplifying to collect all $T(z)$ terms.  Now, Lemma \ref{orthog relation} implies the final sum is nonzero only when $m=-j-k-\ell$.  Performing this substitution, together with collecting all $T(x)$ terms gives
\begin{align*}
D  &=  \frac{1}{(p-1)^3}\sum_{y \in\fpc}\sum_{j,k,\ell=0}^{p-2}G_{-j}G_{-k}G_{-\ell}G_{j+k+\ell}T^j(y^2)T^k(-4)\\
&\qquad \qquad \qquad \cdot T^{-j-k}\left(\frac{27}{1-t}\right)\sum_{x\in\fpc}T^{3k+\ell}(x)\\
  &= \frac{1}{(p-1)^2}\sum_{j,k=0}^{p-2}G_{-j}G_{-k}G_{3k}G_{j-2k}T^k(-4)T^{-j-k}\left(\frac{27}{1-t}\right)\sum_{y\in\fpc}T^{2j}(y).
\end{align*}
The second equality follows by applying the substitution $\ell=-3k$, according to Lemma \ref{orthog relation}, and collecting all $T(y)$ terms. Finally, note that $T^{2j}=\varepsilon$ precisely when $j=0, \frac{p-1}{2}$.  Accounting for both of these cases, we arrive at
\begin{align*} 
D  &= \frac{1}{p-1}\sum_{k=0}^{p-2}G_{0}G_{-k}G_{3k}G_{-2k}T^k(-4)T^{-k}\left(\frac{27}{1-t}\right)\\
  & \qquad \qquad + \frac{1}{p-1}\sum_{k=0}^{p-2}G_{-\frac{p-1}{2}}G_{-k}G_{3k}G_{\frac{p-1}{2}-2k}T^k(-4)T^{-k-\frac{p-1}{2}}\left(\frac{27}{1-t}\right)\\
  &=  \frac{1}{p-1}\sum_{k=0}^{p-2}G_{-k}G_{3k}T^k(-4)\left[-G_{-2k}T^{-k}\left(\frac{27}{1-t}\right)+ \sqrt{p}\,G_{6s-2k}T^{-k-6s}\left(\frac{27}{1-t}\right)\right]\\
  &=  \frac{1}{p-1}\sum_{k=0}^{p-2}G_{-k}G_{3k}T^k(-4)T^{-k}\left(\frac{27}{1-t}\right)\left[-G_{-2k}+\sqrt{p}\,G_{6s-2k}\phi\left(\frac{3}{1-t}\right)\right],
\end{align*}
after collecting like terms and simplifying.  

By a similar analysis, one can compute that
%%%%%%%%%%%%%%%%%%%%%
%A,B, and C Formulas%
%%%%%%%%%%%%%%%%%%%%%
\begin{align*}
A &= -1\\
B &= 1+p\phi\left(\frac{3}{1-t}\right)\\
C &= \frac{1}{p-1}\sum_{j=0}^{p-2}G_{-j}G_{3j}G_{-2j}T^j(-4)T^{-j}\left(\frac{27}{1-t}\right).
\end{align*}
Combining our calculations for $A$, $B$, $C$, and $D$, we see that
\begin{align*}
p\cdot(\#E_t(\fp)-1)&=p^2+A+B+C+D\\
                    &=p^2+p\phi\left(\frac{3}{1-t}\right)\\
& \,\,\qquad +\frac{\sqrt{p}}{p-1}\phi\left(\frac{3}{1-t}\right)\sum_{k=0}^{p-2}G_{-k}G_{3k}G_{6s-2k}T^k(-4)T^{-k}\left(\frac{27}{1-t}\right).
\end{align*}

Now we compute the trace of Frobenius $a(t,p)$. Since $a(t,p)=p+1-\#E_t(\fp)$, we have proved:

%%%%%%%%%%%%%%%%%%%%%%%%%%%%%%%%%
%Statement of formula for a(t,p)%
%%%%%%%%%%%%%%%%%%%%%%%%%%%%%%%%%

\begin{prop}\label{trace calculation}If $p$ is a prime, $p\equiv 1 \pmod{12}$, $s=\frac{p-1}{12}$, and $E_t$ is as in \eqref{defn of E}, then
$$a(t,p)=-\phi\left(\frac{3}{1-t}\right)-\frac{\phi\left(\frac{3}{1-t}\right)}{\sqrt{p}\,(p-1)}\sum_{k=0}^{p-2}G_{-k}G_{3k}G_{6s-2k}T^k(-4)T^{-k}\left(\frac{27}{1-t}\right).$$
\end{prop}

Now that we have a formula for the trace of Frobenius on $E_t$ in terms of Gauss sums, we write our specialization of the $_2F_1$ hypergeometric function in similar terms.  Recall that $s=\frac{p-1}{12}$ and $T$ generates the character group $\widehat{\mathbb{F}_{p}^{\times}}$.  Thus, we may take the character $\xi$ of order 12 in the statement of Theorem \ref{first main theorem} to be $T^s$. %not sure how to explain this%  

The next result gives an explicit formula for $\hg{\xi}{\xi^5}{\varepsilon}{t}$ in terms of Gauss sums.  In its proof, we make use of the specific instances of the Hasse-Davenport relation that were given in Section \ref{preliminaries}.

%%%%%%%%%%%%%%%%%%%%%%
% Computation of 2F1 %
%%%%%%%%%%%%%%%%%%%%%%

\begin{prop}\label{2f1 formula}For $t\in\fp\backslash\{0,1\},$
$$\hg{\xi}{\xi^5}{\varepsilon}{t}=\frac{T^{3s}(4(1-t))}{\sqrt{p}(p-1)} \sum_{k=0}^{p-2}G_{6s-2k}G_{3k}\frac{1}{G_{k}}T^k(4)T^{-k}\left(\frac{27}{1-t}\right).$$
\end{prop}

\begin{proof} By Theorem \ref{greene 44},
%{\allowdisplaybreaks
\begin{align*}
\hg{\xi}{\xi^5}{\varepsilon}{t} &= \xi(-1)\hg{\xi}{\xi^5}{\xi^6}{1-t} & \\
                                    &= T^s(-1)\frac{p}{p-1}\sum_{\chi}\binom{\xi \chi}{\chi}\binom{\xi^5 \chi}{\xi^6 \chi}\chi(1-t) & \text{(Definition \ref{hg})}\\
                                    &=  T^s(-1)\frac{p}{p-1}\sum_{k=0}^{p-2}\binom{T^{s+k}}{T^k}\binom{T^{5s+k}}{T^{6s+k}}T^k(1-t), & 
\end{align*}
%}
as $T$ generates the character group.  Now we rewrite the product $\binom{T^{s+k}}{T^k}\binom{T^{5s+k}}{T^{6s+k}}$ of binomial coefficients in terms of Gauss sums, by way of Lemma \ref{greene 29}.  Since $T^s=\xi$ and $T^{-s}=\xi^{-1}$ are not trivial, we have

%{\allowdisplaybreaks
\begin{align*}
\binom{T^{s+k}}{T^k}\binom{T^{5s+k}}{T^{6s+k}} &= \left[\frac{G_{s+k}G_{-k}T^k(-1)}{pG_{s}}\right] \cdot \left[ \frac{G_{5s+k}G_{-6s-k}T^{6s+k}(-1)}{pG_{-s}} \right]\\
                                               &= \frac{1}{p^3}G_{s+k}G_{-k}G_{5s+k}G_{-6s-k}T^{5s+2k}(-1),
\end{align*}
%}
since $G_sG_{-s}=pT^s(-1)$ by Lemma \ref{greene 112}.  Thus, 
%{\allowdisplaybreaks
\begin{align*}
\hg{\xi}{\xi^5}{\varepsilon}{t} &= \frac{T^s(-1)}{p^2(p-1)} \sum_{k=0}^{p-2} G_{s+k}G_{-k}G_{5s+k}G_{-6s-k}T^{5s+2k}(-1)T^k(1-t)\\
                                &= \frac{\phi(-1)}{p^2(p-1)} \sum_{k=0}^{p-2} G_{s+k}G_{-k}G_{5s+k}G_{-6s-k}T^k(1-t),
\end{align*}
%}
since $T^sT^{5s}=\phi$ and $T^{2k}(-1)=1 \,\,\mbox{for all}\, k.$

Now we apply the Hasse-Davenport relation (Corollary \ref{HD2}) and make a substitution for $G_{-k}G_{-6s-k}$.  We obtain
$$\hg{\xi}{\xi^5}{\varepsilon}{t}=\frac{\phi(-1)}{p^{\frac{3}{2}}(p-1)} \sum_{k=0}^{p-2} G_{s+k}G_{5s+k}G_{-2k}T^k(4)T^k(1-t).$$

\noindent Next, we let $k \mapsto k+3s$ and find  
%{\allowdisplaybreaks
\begin{align*}
\hg{\xi}{\xi^5}{\varepsilon}{t} &= \frac{\phi(-1)}{p^{\frac{3}{2}}(p-1)} \sum_{k=0}^{p-2} G_{4s+k}G_{8s+k}G_{-2k-6s}T^{k+3s}(4)T^{k+3s}(1-t)\\
                                &= \frac{\phi(-1)T^{4s}(-1)}{\sqrt{p}(p-1)} \sum_{k=0}^{p-2} G_{6s-2k}G_{3k}\frac{1}{G_k}T^{-k}(27)T^{k+3s}(4)T^{k+3s}(1-t),
\end{align*}
%}
by applying the Hasse-Davenport relation (Corollary \ref{HD3}) to make a substitution for $G_{4s+k}G_{8s+k}$, and by noting that $G_{-2k-6s}=G_{-2k+6s}$.  Then, since  $p \equiv 1\pmod{12}$ implies $\phi(-1)T^{4s}(-1)=T^{10s}(-1)=1$, we simplify to obtain
\begin{equation*}
\hg{\xi}{\xi^5}{\varepsilon}{t} = \frac{T^{3s}(4(1-t))}{\sqrt{p}(p-1)} \sum_{k=0}^{p-2}G_{6s-2k}G_{3k}\frac{1}{G_{k}}T^k(4)T^{-k}\left(\frac{27}{1-t}\right),
\end{equation*}
as desired.
\end{proof}

%%%%%%%%%%%%%%%%%%%%%
% proof of theorem%%%
%%%%%%%%%%%%%%%%%%%%%

We now have the necessary tools to complete the proof of Theorem \ref{first main theorem}.  

\begin{proof}[Proof of Theorem \ref{first main theorem}]
We combine the results of Propositions \ref{trace calculation} and \ref{2f1 formula} with a bit of algebra to complete the proof.  We begin by taking the formula for $\hg{\xi}{\xi^5}{\varepsilon}{t}$ given in Proposition \ref{2f1 formula}, splitting off the $k=0$ term in the sum and applying Lemma \ref{greene 112} to the $k\geq1$ terms, to move all Gauss sums to the numerator.  We also simplify by noticing that  $T^{-k}(-1)=T^k(-1)$ implies $T^{-k}(-1)T^k(4)=T^k(-4)$.  We see that
\begin{align*}
\hg{\xi}{\xi^5}{\varepsilon}{t} &= \frac{T^{3s}(4(1-t))}{\sqrt{p}(p-1)} \Biggl[ \sqrt{p}+  \frac{1}{p}\sum_{k=1}^{p-2} G_{6s-2k}G_{3k}G_{-k}T^{-k}(-1)T^k(4)\Biggr.\\
&\Biggl. \hspace*{0.4in} \cdot T^{-k}\left( \frac{27}{1-t} \right) \Biggr]\\
                                &= \frac{T^{3s}(4(1-t))}{\sqrt{p}(p-1)} \left[ \sqrt{p} + \frac{1}{p}\sum_{k=1}^{p-2} G_{6s-2k}G_{3k}G_{-k}T^k(-4)T^{-k}\left( \frac{27}{1-t} \right) \right].
\end{align*}
Next, we multiply  by $\displaystyle{\frac{\phi(3)T^{3s}(1-t)}{\phi(3)T^{3s}(1-t)}}$ and rearrange, while recalling that $\phi=\phi^{-1}$.  We obtain
\begin{align*}
\hg{\xi}{\xi^5}{\varepsilon}{t} &= -\frac{T^{3s}(4)}{\phi(3)T^{3s}(1-t)} \left[ -\frac{1}{p-1}\phi \left( \frac{3}{1-t} \right)\right.\\ 
                                & \hspace*{0.2in}\left.- \frac{1}{p^{\frac{3}{2}}(p-1)}\phi \left( \frac{3}{1-t} \right) \sum_{k=1}^{p-2} G_{6s-2k}G_{3k}G_{-k}T^k(-4)T^{-k}\left( \frac{27}{1-t} \right) \right]\\
                                &= -T^{3s}(4)\phi(3)T^{-3s}(1-t) \left[ -\frac{\phi \left( \frac{3}{1-t} \right)}{p} \right.\\
                                & \hspace*{0.2in}\left. -  \frac{\phi \left( \frac{3}{1-t} \right)}{p^{\frac{3}{2}}(p-1)}  \sum_{k=0}^{p-2} G_{6s-2k}G_{3k}G_{-k}T^k(-4)T^{-k}\left( \frac{27}{1-t} \right) \right].
\end{align*}
\noindent The last equality follows by noting that the $k=0$ term of the final sum is $-\frac{\phi \left( \frac{3}{1-t} \right)}{p(p-1)}$ and 
$$ -\frac{\phi \left( \frac{3}{1-t} \right)}{p-1}+ \frac{\phi \left( \frac{3}{1-t} \right)}{p(p-1)}= -\frac{\phi \left( \frac{3}{1-t} \right)}{p}.$$

%so
%$$ -\frac{\phi \left( \frac{3}{1-t} \right)}{p} - \frac{\phi \left( \frac{3}{1-t} \right)}{p(p-1)} = -\frac{\phi \left( \frac{3}{1-t} \right)}{p-1}.$$
\noindent  Recalling the expression for $a(t,p)$ given by Proposition \ref{trace calculation}, we have that

%Recall now that Proposition \ref{trace calculation} provides that 
%$$-\phi\left(\frac{3}{1-t}\right)-\frac{\phi\left(\frac{3}{1-t}\right)}{\sqrt{p}(p-1)}\sum_{k=0}^{p-2}G_{6s-2k}G_{3k}G_{-k}T^k(-4)T^{-k}\left(\frac{27}{1-t}\right)=a(t,p).$$ Thus, we have 
$$p\,\hg{\xi}{\xi^5}{\varepsilon}{t}=-T^{3s}(4)\phi(3)T^{-3s}(1-t)a(t,p),$$
so the proof is complete if $T^{3s}(4)\phi(3)T^{-3s}(1-t)=\phi(2)\xi^{-3}(1-t)$.  Since $T^{3s}=\xi^3$ and  $T^{-3s}=\xi^{-3}$, we need only show that 
\begin{equation}\label{simplify character}
\xi^3(4)\phi(3)=\phi(2).
\end{equation}
By multiplicativity, $\xi^3(4)=\xi^6(2)=\phi(2)$.  Further, $\phi(3)=\bigl(\frac{3}{p}\bigr)=\bigl(\frac{p}{3}\bigr)$ by quadratic reciprocity, since $p\equiv1\pmod{4}$.  Also, since $p\equiv1\pmod{3}$, we have $\phi(3)=\left(\frac{1}{3}\right)=1$.  This verifies (\ref{simplify character}), and hence completes the proof.
\end{proof}

%%%%%%%%%%%%%%%%%%%
%Other two results%
%%%%%%%%%%%%%%%%%%%

We have proved two other results similar to Theorem \ref{first main theorem}, but which apply to different families of elliptic curves.  These results are finite field analogues of Beukers results (see \cite{Be93}) relating periods of families of elliptic curves to values of classical hypergeometric functions, as described in Section \ref{history hg and ec}.  Note that the characters which appear in our $_2F_1$ bear a striking resemblance to the parameters Beukers used in the classical case.

%The particular families of elliptic curves were considered by Beukers in \cite{Be93}, where he related the periods of these families to values of classical hypergeometric functions.  Our results are analogues involving hypergeometric functions over $\fp$, and the characters which appear in our $_2F_1$ bear a striking resemblance to the parameters Beukers used in the classical case.  We now state these results without proof, as they are proved following the same general steps given in the proof of Theorem \ref{first main theorem}.\\

\begin{prop}\label{prop 2.5}
Suppose $p \equiv 1 \pmod{12}$ is prime, and let  $\xi\in \widehat{\mathbb{F}_{p}^{\times}}$ have order $12$.  Let $E_t:y^2=x^3+tx+1$, and let $a(t,p)=p+1-\#E_t(\fp)$.  Then
$$p \, \hg{\xi}{\xi^7}{\xi^8}{-\frac{4}{27}t^3}=\chi(t)a(t,p),$$
where $\chi(t)=-\xi^{-1}(-4)\xi^{-4}(\frac{t^3}{27})$.
\end{prop}

\begin{prop}\label{prop 2.4}
Suppose $p \equiv 1 \pmod{12}$ is prime, and let  $\xi\in \widehat{\mathbb{F}_{p}^{\times}}$ have order $12$.  Let $E_t:y^2=x^3-x-t$, and let $a(t,p)=p+1-\#E_t(\fp)$.  Then
$$p \, \hg{\xi}{\xi^5}{\phi}{\frac{27}{4}t^2}=-\xi^3(-27)a(t,p).$$
\end{prop}

It is interesting to note that in Proposition \ref{prop 2.5}, the values of the character $\chi(t)$, which appears as the coefficient of $a(t,p)$, are simply sixth roots of unity, and in Proposition \ref{prop 2.4}, the values of $\xi^3(-27)$ are simply $\pm 1$.  A priori, $\xi^3(-27)\in\{\pm1,\pm i\},$ but in fact, we have $(\xi^3(-27))^2=(\xi^3(-1)\xi^3(27))^2=\phi(-1)\phi(27)=\phi(-1)\phi(3)=1$.  This follows since $p\equiv 1\pmod{12}$ implies $\phi(-1)=1$ and since $\phi(3)=1$, as shown in the proof of Theorem \ref{first main theorem}.

%%%%%%%%%%%%%%%%%%%%%%%%%%
%% histosry of HGF and MF%
%%%%%%%%%%%%%%%%%%%%%%%%%%

\section{Recent History: Hypergeometric Functions and Modular Forms}\label{history hg and mf}

As with elliptic curves, classical hypergeometric functions have connections to modular forms.  Investigations into these relations began in the early twentieth century.  More recently, Stiller \cite{St88} proved an array of results linking the two objects.  In fact, the classical hypergeometric series and family of elliptic curves (i.e. one with $j$-invariant $\frac{1728}{t}$) considered by Stiller prompted the choice of family $E_t$ and the $_2F_1$ function used in our main results.  For this reason, we now state one of Stiller's results in full.  We let 

$$E_4(q)=1+204\sum_{n\geq 1} \sigma_3(n)q^n$$
and
$$E_6(q)=1-504 \sum_{n\geq 1} \sigma_5(n)q^n$$
be the classical Eisenstein series of weights $4$ and $6$, respectively, for $\Gamma$.  Stiller directly related these two modular forms to classical hypergeometric series:

\begin{theorem}[Stiller \cite{St88}, Thm. 5]
Let $\C[E_4,E_6]$ be the graded algebra of modular forms for $\Gamma$ and let $\C [ _2F_1[\frac{1}{12},\frac{5}{12};1;t]^4, (1-t)^{1/2} \, _2F_1[\frac{1}{12},\frac{5}{12};1;t]^6]$ be the graded algebra of hypergeometric functions (graded by the power of $_2F_1$).  Then these two algebras are canonically isomorphic as graded algebras of power series in $q=e^{2\pi i z}$ and $t$, respectively.  Moreover, the isomorphism is the pull-back $\pi^*$, where $\pi(q)=\frac{1728}{j(q)}$ and $j$ is the usual elliptic modular function. 
\end{theorem}

%NEED TO FINISH THIS OFF AND MAKE IT MAKE SENSE FOR OUR CONTEXT.  Need to:

%\begin{itemize}
%\item Describe the pull back so the theorem statement is self contained in this paper.
%\item Need to include some information about the family of elliptic curves.  Specifically, need to find where Stiller considers a family of j invariant 1/t in this paper....having trouble with this!
%\end{itemize}

Following Greene's introduction of hypergeometric functions over finite fields, results emerged relating them to modular forms.  Ahlgren, Ono, and others produced formulas for traces of Hecke operators on certain spaces of cusp forms.  These formulas were given in terms of traces of Frobenius on related families of elliptic curves.  %However, Theorem \ref{koike ono} implies further relationships to hypergeometric functions.

Specifically, in 2000 and 2002, Ahlgren and Ono \cite{AO00} and Ahlgren \cite{Ah02} exhibited formulas for the traces of Hecke operators on spaces of cusp forms in levels 8 and 4.  Let $k\geq2$ be an even integer, and define
$$F_k(x,y)=\frac{x^{k-1}-y^{k-1}}{x-y}.$$
Then letting $x+y=s$ and $xy=p$ gives rise to polynomials $G_k(s,p)=F_k(x,y)$.  These polynomials can be written alternatively as
\begin{equation}\label{defn of G}
G_k(s,p)=\sum_{j=0}^{\frac{k}{2}-1} (-1)^j \binom{k-2-j}{j} p^j s^{k-2j-2}.
\end{equation}
The results given below are given in \cite{AO00} and \cite{Ah02} for the cases of level 4, weight 6 and of level 8, weight 4, respectively.  However, the statements hold for all even $k\geq 4$ with the same proofs.

\begin{theorem}[(a) Ahlgren and Ono \cite{AO00}, (b) Ahlgren \cite{Ah02}]
Let $p$ be an odd prime and $k\geq4$ be an even integer.  Then

$(a) \quad \displaystyle{\,\,\textnormal{Tr}_k(\Gamma_0(8),p)=-4-\sum_{t=2}^{p-2} G_k(_2A_1(p,t^2),p)}$
	
$(b) \quad \displaystyle{\,\,\textnormal{Tr}_k(\Gamma_0(4),p)=-3-\sum_{t=2}^{p-1} G_k(_2A_1(p,t),p).}$
\end{theorem}

Ahlgren and Ono's methods involved combining the Eichler-Selberg trace formula \cite{Hi89} with a theorem given by Schoof \cite{Sc87}.  In the proof of Theorem \ref{second main theorem}, given in the next section, we use similar techniques to exhibit a formula in the level 1 setting.  Recently, Frechette, Ono, and Papanikolas expanded the techniques of Ahlgren and Ono and obtained results in the level 2 case:

\begin{theorem}[Frechette, Ono, and Papanikolas \cite{FOP04}]
Let $p$ be an odd prime and $k\geq 4$ be even.  When $p \equiv 1 \pmod{4}$, write $p=a^2+b^2$, where $a,b$ are nonnegative integers, with $a$ odd.  Then
$$\textnormal{Tr}_k(\Gamma_0(2),p)=-2-\delta_k(p)-\sum_{t=1}^{p-2} G_k(_3A_2(p,t),p),$$
where 
\begin{equation*}
\delta_k(p)=
\begin{cases}
\frac{1}{2}G_k(2a,p)+\frac{1}{2}G_k(2b,p) &\textnormal{if}\,\, p\equiv1\pmod{4}\\
(-p)^{k/2-1} &\textnormal{if}\,\, p\equiv3\pmod{4}.
\end{cases}
\end{equation*}
\end{theorem}

In addition, Frechette, Ono, and Papanikolas used relationships between counting points on varieties over $\fp$ and hypergeometric functions over $\fp$ to obtain further results for the traces of Hecke operators on spaces of newforms in level 8.  Most recently, Papanikolas \cite{Pa06} used the results in \cite{FOP04} as a starting point to obtain a new formula for Ramanujan's $\tau$ function, as well as a new congruence for $\tau(p) \pmod{11}$.
%(STATE THIS??)

\section{Proof of Theorem \ref{second main theorem}}\label{proof of theorem 2}

%%%%%%%%%%%%%%%%%%%%%%%%%%%%%%%%%%%%%%%%%%
%% PRELIMINARIES FOR PROOF OF THEOREM 2 %%
%%%%%%%%%%%%%%%%%%%%%%%%%%%%%%%%%%%%%%%%%%

The proof of Theorem \ref{second main theorem} utilizes three important results.  First, we use Hasse's classical bound on the number of points on an elliptic curve defined over a finite field.  (See, for example, \cite{Si86} page 131 for details.)  We also use a theorem of Schoof, together with Hijikata's version of the Eichler-Selberg trace formula, which require some notation.  We follow the treatment given in \cite{FOP04}.  If $d<0$, $d\equiv 0,1 \pmod{4}$, let $\mathcal{O}(d)$ denote the unique imaginary quadratic order in $\mathbb{Q}(\sqrt{d})$ having discriminant $d$.  Let $h(d)=h(\mathcal{O}(d))$ be the order of the class group of $\mathcal{O}(d)$, and let $w(d)=w(\mathcal{O}(d))$ be half the cardinality of the unit group of $\mathcal{O}(d)$.  We then let $h^*(d)=h(d)/w(d)$.  Further, if $d$ is the discriminant of an imaginary quadratic order $\mathcal{O}$, let
\begin{equation}\label{defn of H}
H(d):=\sum_{\mathcal{O}\subseteq\mathcal{O}'\subseteq\mathcal{O}_{max}}h(\mathcal{O}'),
\end{equation}
where the sum is over all orders $\mathcal{O}'$ between $\mathcal{O}$ and $\mathcal{O}_{max}$, the maximal order.  A complete treatment of the theory of orders in imaginary quadratic fields can be found in section 7 of \cite{Co89}.

Additionally, if $K$ is a field, we define
$$\emph{Ell}_K:=\{[E]_K | E \,\textnormal{is defined over}\, K\},$$
where $[E]_K$ denotes the isomorphism class of $E$ over $K$ and $[E_1]_K=[E_2]_K$ if there exists an isomorphism $\beta:E_1 \rightarrow E_2$ over $K$.  Now if $p$ is an odd prime, define 
\begin{equation}\label{defn of I}
I(s,p):=\{[E]_{\fp}\in\emph{Ell}_{\fp} | \#E(\fp)=p+1 \pm s\}.
\end{equation}
Schoof proved the following theorem, connecting the quantities in (\ref{defn of H}) and (\ref{defn of I}).

\begin{theorem}[Schoof \cite{Sc87}, Thm. 4.6]\label{schoof's theorem}  If $p$ is an odd prime and $s$ is an integer with $0<s<2\sqrt{p}$, then
\[
\#I(s,p)=2H(s^2-4p).
\]
\end{theorem}

The final key ingredient to the proof of Theorem \ref{second main theorem} is the Eichler-Selberg trace formula, which provides a starting point for calculating the trace of the $p^{th}$ Hecke operator on $S_k$.  We use Hijikata's version of this formula, which is found in \cite{Hi89}, but we only require the level 1 formulation.  

Let $p\equiv 1 \pmod{12}$ be prime, and recall the definition of polynomials $G_k(s,p)$ given in (\ref{defn of G}), and in the statement of Theorem \ref{second main theorem}.  Note that when writing the polynomials $G_k(s,p)$, we take the convention $s^0=1$, so that the constant term of $G_k(s,p)$ is $(-p)^{\frac{k}{2}-1}$, for all values of $s$.  Using this notation, the formulation given below is a straightforward reduction of Hijikata's trace formula in the level one case.

\begin{theorem}[Hijikata \cite{Hi89}, Thm. 2.2]\label{hijikata}
Let $k\geq 2$ be an even integer, and let $p\equiv 1 \pmod{12}$ be prime.  Then
$$\textnormal{Tr}_k(\Gamma,p)=-h^*(-4p)(-p)^{\frac{k}{2}-1}-1-\sum_{0<s<2\sqrt{p}}G_k(s,p)\sum_{f|\ell}h^*\left(\frac{s^2-4p}{f^2}\right)+\delta(k),$$
where 
$$\delta(k)=\begin{cases} p+1 &\mbox{if}\,\, k=2\\0 & \mbox{otherwise}\end{cases}$$ 
and where we classify integers $s$ with $s^2-4p<0$ by some positive integer $\ell$ and square-free integer $m$ via 
$$s^2-4p=\begin{cases}\ell^2m, & 0>m\equiv 1 \pmod{4}\\\ell^24m, & 0>m\equiv 2,3 \pmod{4}.\end{cases}$$

\end{theorem}

Next, we recall a result which relates isomorphism classes in $\emph{Ell}_{\overline{\mathbb{F}}_p}$ and $\emph{Ell}_{\fp}$.  Define a map
\begin{align*}
\eta:\emph{Ell}_{\fp}&\rightarrow \emph{Ell}_{\overline{\mathbb{F}}_p}\\
[E]_{\fp} &\mapsto [E]_{\overline{\mathbb{F}}_p}.
\end{align*}
 Note that $\eta$ is well defined since two curves which are isomorphic over $\fp$ are necessarily isomorphic over $\overline{\mathbb{F}}_p$.

\begin{lemma}[]\label{fp and fpbar}
Let $p\geq5$ be prime.  Suppose $[E]_{\overline{\mathbb{F}}_p}\in\emph{Ell}_{\overline{\mathbb{F}}_p}$ and $E$ is defined over $\fp$.  Then
$$\#\eta^{-1}([E]_{\overline{\mathbb{F}}_p})=\begin{cases} 2 & \textnormal{if} \, j\neq 0, 1728\\ 4 & \textnormal{if} \, j=1728 \\ 6 & \textnormal{if} \, j=0.\end{cases}$$
\end{lemma}
\begin{proof}
See Section X.5 of \cite{Si86}.
\end{proof}

Among isomorphism classes of elliptic curves over $\fp$, two are of particular interest to us: those having $j$-invariant $1728$ and those having $j$-invariant $0$.  If $E$ is any elliptic curve defined over $\fp$, we let $a(E)$ be given by $a(E)=p+1-\#E(\fp)$. The following two lemmas compute formulas for the sums of $a(E)^n$ over all curves $E$ over $\fp$ having $j$-invariant $1728$ or $0$, respectively.

\begin{lemma}[\cite{Fu07}, Lemma IV.3.3]\label{j inv 1728}
Let $p\equiv 1 \pmod{12}$ and let $a,b\in\mathbb{Z}$ be such that  $p=a^2+b^2$ and $a+bi \equiv 1 \,(2+2i)$ in $\mathbb{Z}[i]$.  Then for $n\geq2$ even,
$$\sum_{\substack{[E]_{\fp}\in Ell_{\fp}\\j(E)=1728}}a(E)^n=2^{n+1}(a^n+b^n).$$
\end{lemma}

\begin{lemma}[\cite{Fu07}, Lemma IV.3.5]\label{j inv 0}
Let $p\equiv 1 \pmod{12}$ and let $c,d\in\mathbb{Z}$ such that $p=c^2-cd+d^2$ and $c+d\omega \equiv 2 \,(3)$ in $\mathbb{Z}[\omega]$, where $\omega=e^{2\pi i /3}$.  Then for $n\geq2$ even,
$$\sum_{\substack{[E]_{\fp}\in Ell_{\fp}\\j(E)=0}}a(E)^n=2[(c+d)^n+(2c-d)^n+(c-2d)^n].$$
\end{lemma}

We omit the proofs of Lemmas \ref{j inv 1728} and \ref{j inv 0}, as they are quite tedious and require checking dozens of cases.  However, these proofs are not difficult.  The only tools used are the classification of elliptic curves over $\fp$ having $j$-invariant 1728 or 0, together with known formulas for $a(E)$ in these cases.  These formulas (see Chapter 18 of \cite{IR90}) are given in terms of $m^{th}$ power residues, whose values must be calculated on a case by case basis for each curve with the given $j$-invariant.  Complete details of the proofs of Lemmas \ref{j inv 1728} and \ref{j inv 0} can be found in \cite{Fu07}.

Now we proceed toward the proof of Theorem \ref{second main theorem}.  As before, let $p\equiv 1\pmod{12}$ be prime.  As in (\ref{defn of E}), we define a family of elliptic curves over $\fp$ by
$$E_t:y^2=4x^3-\frac{27}{1-t}x-\frac{27}{1-t}.$$  Further, for $t\in\fp$, $t\neq 0,1$, recall that
$$a(t,p)=p+1-\#E_t(\fp).$$
As in  Lemmas \ref{j inv 1728} and \ref{j inv 0}, we let integers $a,b,c,$ and $d$ be defined by $p=a^2+b^2=c^2-cd+d^2$, where $a+bi\equiv 1 \, (2+2i)$ in $\mathbb{Z}[i]$ and  $c+d\omega\equiv 2 \, (3)$ in $\mathbb{Z}[\omega]$, where $\omega=e^{2\pi i/3}$.  Finally, we let  $h$, $h^{*}$, $w$, and $H$ be defined as at the start of this section.  

\begin{lemma}\label{h to h*} If $p\equiv 1 \pmod{12}$ is prime and notation is as above, then for $n\geq 2$ even,
\begin{multline*}
\sum_{0<s<2\sqrt{p}}s^n \sum_{f|\ell}h\left(\frac{s^2-4p}{f^2}\right)=\sum_{0<s<2\sqrt{p}}s^n \sum_{f|\ell}h^*\left(\frac{s^2-4p}{f^2}\right)\\+\frac{1}{4}\sum_{\substack{[E]_{\fp}\in \emph{Ell}_{\fp}\\j(E)=1728}}a(E)^n+\frac{1}{3}\sum_{\substack{[E]_{\fp}\in \emph{Ell}_{\fp}\\j(E)=0}}a(E)^n,
\end{multline*}
where we classify integers $s$ with $s^2-4p<0$ by some positive integer $\ell$ and square-free integer $m$ via 
$$s^2-4p=\begin{cases}\ell^2m, & 0>m\equiv 1 \pmod{4}\\\ell^24m, & 0>m\equiv 2,3 \pmod{4}.\end{cases}$$
\end{lemma}

\begin{proof}
First, notice that $h$ and $h^*$ agree unless the argument $\frac{s^2-4p}{f^2}=-3$ or $-4$, since in all other cases $w(d)=1$.  Thus, we have
\begin{multline*}
\sum_{0<s<2\sqrt{p}}s^n \sum_{f|\ell}h\left(\frac{s^2-4p}{f^2}\right)=\sum_{0<s<2\sqrt{p}}s^n\sum_{\substack{f|\ell\\\frac{s^2-4p}{f^2}\neq -3,-4}}h^{*}\left(\frac{s^2-4p}{f^2}\right)\\
+\sum_{0<s<2\sqrt{p}}s^n\sum_{\substack{f|\ell\\\frac{s^2-4p}{f^2}=-4}}h(-4)+\sum_{0<s<2\sqrt{p}}s^n\sum_{\substack{f|\ell\\\frac{s^2-4p}{f^2}=-3}}h(-3).
\end{multline*}
When $\frac{s^2-4p}{f^2}=-4$, we have the maximal order $\mathbb{Z}[i]$ and  $h^*(-4)=\frac{h(-4)}{w(-4)}=\frac{1}{2}$, so $h(-4)=h^*(-4)+\frac{1}{2}$.  On the other hand, when $\frac{s^2-4p}{f^2}=-3$, we have the maximal order $\mathbb{Z}[\omega]$ and $h^*(-3)=\frac{h(-3)}{w(-3)}=\frac{1}{3}$, so $h(-3)=h^*(-3)+\frac{2}{3}$.  Making these substitutions, we see that
\begin{multline}\label{h* step}
\sum_{0<s<2\sqrt{p}}s^n \sum_{f|\ell}h\left(\frac{s^2-4p}{f^2}\right)=\sum_{0<s<2\sqrt{p}}s^n \sum_{f|\ell}h^*\left(\frac{s^2-4p}{f^2}\right)\\
+\frac{1}{2}\sum_{0<s<2\sqrt{p}}s^n\sum_{\substack{f|\ell\\\frac{s^2-4p}{f^2}=-4}}1\,\,+\frac{2}{3}\sum_{0<s<2\sqrt{p}}s^n\sum_{\substack{f|\ell\\\frac{s^2-4p}{f^2}=-3}}1.
\end{multline}
To complete the proof, we must verify that
\begin{equation}\label{h=-4}
\sum_{0<s<2\sqrt{p}}s^n\sum_{\substack{f|\ell\\\frac{s^2-4p}{f^2}=-4}}1 = \frac{1}{2}\sum_{\substack{[E]_{\fp}\in\emph{Ell}_{\fp}\\j(E)=1728}}a(E)^n
\end{equation}
and
\begin{equation}\label{h=-3}
\sum_{0<s<2\sqrt{p}}s^n\sum_{\substack{f|\ell\\\frac{s^2-4p}{f^2}=-3}}1 = \frac{1}{2}\sum_{\substack{[E]_{\fp}\in\emph{Ell}_{\fp}\\j(E)=0}}a(E)^n.
\end{equation}

First, we consider \eqref{h=-4}.  Using known formulas for $a(E)$ in the case of curves with $j$-invariant 1728 (see Chapter 18 of \cite{IR90}), one can show that $a(E)=\pm 2a,\pm 2b$ for all $E$ relevant to \eqref{h=-4} with $j$-invariant $1728$.  Also, it is easy to verify that $s=|2a|,|2b|$ satisfy $\frac{s^2-4p}{\ell^2}=-4$ (with $\ell=|b|,|a|$, respectively).  Now, suppose $0<s<2\sqrt{p}\,$ satisfies  $\frac{s^2-4p}{\ell^2}=-4$.  Then $s^2-4p=-4\ell^2$ implies $s$ is even, so we have $\left(\frac{s}{2}\right)^2+\ell^2=p.$  Thus, it must be that $\frac{s}{2}=|a|,|b|$, since $\mathbb{Z}[i]$ is a UFD and $p=a^2+b^2$.  Since $n$ is even, $(2a)^n=(-2a)^n$ and $(2b)^n=(-2b)^n$, so \eqref{h=-4} follows.

We prove \eqref{h=-3} in a similar manner.  Using known formulas for $a(E)$ in the case of curves with $j$-invariant 0 (see Chapter 18 of \cite{IR90}), one can show that $a(E)=\pm(c+d),\pm(2c-d),\pm(c-2d)$ for all $E$ with $j$-invariant $0$ that appear in \eqref{h=-3}.  Also, $s=|c+d|$, $|2c-d|$, and $|c-2d|$ satisfy  $\frac{s^2-4p}{\ell^2}=-3$ (by taking $\ell=|c-d|$, $|d|$, and $|c|$, respectively).  Now, suppose  $\frac{s^2-4p}{\ell^2}=-3$.  Then in $\mathbb{Z}[\sqrt{-3}]$, we have 
$$4p=(s+\sqrt{-3}\,\ell)(s-\sqrt{-3}\,\ell).$$ 
Since $-3\equiv 5 \pmod{8}$, $2$ is inert in $\mathbb{Z}[\sqrt{-3}]$, so we must have $2|(s\pm \sqrt{-3}\,\ell)$.  This implies
\begin{equation}\label{factor p}
p=\left(\frac{s}{2}+\sqrt{-3}\,\frac{\ell}{2}\right)\left(\frac{s}{2}-\sqrt{-3}\,\frac{\ell}{2}\right)
\end{equation}
in $\mathbb{Z}[\sqrt{-3}]$.  Recall that we have $p=c^2-cd+d^2$.  In $\mathbb{Z}[\omega]$, we can write this as
\begin{align}
p&=c^2-cd+d^2 \notag\\
&=(c+d\omega)(c+d\omega^2) \label{first}\\
&=(d+c\omega)(d+c\omega^2)\label{second}\\
&=(c\omega+d\omega^2)(c\omega^2+d\omega).\label{third}
\end{align}
Since $\omega=e^{2\pi i/3}=-\frac{1}{2}+\frac{\sqrt{-3}}{2}$, we can consider each of these factorizations in $\mathbb{Z}[\sqrt{-3}]$, and each must be the same as \eqref{factor p}, since $\mathbb{Z}[\sqrt{-3}]$ is a UFD.  Making the substitution for $\omega$ into \eqref{first}, \eqref{second}, and \eqref{third} and comparing to \eqref{factor p} implies that $s=|2c-d|$, $|2d-c|$, and $|c+d|$, respectively.  So in fact,  $s=|2c-d|$, $|2d-c|$, $|c+d|$ are the only contributing $s$ values to the sum on the left hand side of \eqref{h=-3}.  Then since $a(E)=\pm(c+d),\pm(2c-d),\pm(c-2d)$ and $n$ is even, we have proved \eqref{h=-3}.

The lemma is finally proved by making the substitutions from \eqref{h=-4} and \eqref{h=-3} into \eqref{h* step}. 
\end{proof}

\begin{prop}\label{step 1 of thm 2 proof}
Let $p\equiv 1 \pmod{12}$ be prime and notation as above.  Then for $n\geq2$ even, 
\begin{multline*}\sum_{t=2}^{p-1}a(t,p)^n=\sum_{0<s<2\sqrt{p}}s^n \sum_{f|\ell}h^{*} \left(\frac{s^2-4p}{f^2}\right)\\-2^{n-1}(a^n+b^n)-\frac{1}{3}[(c+d)^n+(2c-d)^n+(c-2d)^n],\end{multline*}where we classify integers $s$ with $s^2-4p<0$ by some positive integer $\ell$ and square-free integer $m$ via 
$$s^2-4p=\begin{cases}\ell^2m, & 0>m\equiv 1 \pmod{4}\\\ell^24m, & 0>m\equiv 2,3 \pmod{4}.\end{cases}$$
\end{prop}

\begin{proof}
Notice that for the given family of elliptic curves, $j(E_t)=\frac{1728}{t}$.  Thus, as $t$ ranges from $2$ to $p-1$, each $E_t$ represents a distinct isomorphism class of elliptic curves in $\emph{Ell}_{\overline{\mathbb{F}}_p}$.  Moreover, since $j(E_t)$ gives an automorphism of $\mathbb{P}^1$, every $j$-invariant other than $0$ and $1728$ is represented precisely once.  Thus, for even $n\geq 2$, we have
$$\sum_{t=2}^{p-1}a(t,p)^n=\sum_{\substack{[E]_{\overline{\mathbb{F}}_p}\in Ell_{\overline{\mathbb{F}}_p}\\E/\fp\\j(E)\neq 0,1728}}a(E)^n.$$

For elliptic curves with $j$-invariant other than $0$ and $1728$, each class $[E]\in\emph{Ell}_{\overline{\mathbb{F}}_p}$ gives rise to two distinct classes in $\emph{Ell}_{\fp}$ (see Lemma \ref{fp and fpbar}), represented by $E$ and its quadratic twist $E^{tw}$.  For such curves, $a(E)$ and $a(E^{tw})$ differ only by a sign, and so $a(E)^n=a(E^{tw})^n$, since $n$ is even.  Therefore, we have
$$\sum_{t=2}^{p-1}a(t,p)^n=\sum_{\substack{[E]_{\overline{\mathbb{F}}_p}\in Ell_{\overline{\mathbb{F}}_p}\\E/\fp\\j(E)\neq 0,1728}}a(E)^n=\frac{1}{2}\sum_{\substack{[E]_{\fp}\in Ell_{\fp}\\j(E)\neq 0,1728}}a(E)^n.$$
Then, if we add and subtract the contributions from the classes $[E]_{\fp}\in\emph{Ell}_{\fp}$ with $j(E)=0, 1728$, we have
\begin{equation}\label{prop step 2}
\sum_{t=2}^{p-1}a(t,p)^n=\frac{1}{2}\left[ \sum_{[E]_{\fp}\in Ell_{\fp}}a(E)^n - \sum_{\substack{[E]_{\fp}\in Ell_{\fp}\\j(E)=1728}}a(E)^n-\sum_{\substack{[E]_{\fp}\in Ell_{\fp}\\j(E)=0}}a(E)^n\right].
\end{equation}

Now we look more closely at the sum $\displaystyle{\sum_{[E]_{\fp}\in \emph{Ell}_{\fp}}a(E)^n}$.  By Hasse's theorem, $\emph{Ell}_{\fp}$ is the the disjoint union 
$$\emph{Ell}_{\fp}=\bigcup_{0\leq s<2\sqrt{p}}I(s,p),$$
where $I(s,p)$ is defined as in (\ref{defn of I}).  Then since $n\geq2$ is even, we may write
\begin{align*}
\sum_{[E]_{\fp}\in \emph{Ell}_{\fp}}a(E)^n&=\sum_{0\leq s <2\sqrt{p}}\quad\sum_{[E]_{\fp}\in I(s,p)} s^n\\
&=\sum_{0< s <2\sqrt{p}}\#I(s,p)s^n,
\end{align*}
since $s=0$ makes no contribution.  Substituting this into (\ref{prop step 2}) gives
$$\sum_{t=2}^{p-1}a(t,p)^n=\frac{1}{2}\sum_{0<s<2\sqrt{p}}\#I(s,p)s^n-\frac{1}{2}\sum_{\substack{[E]_{\fp}\in \emph{Ell}_{\fp}\\j(E)=1728}}a(E)^n-\frac{1}{2}\sum_{\substack{[E]_{\fp}\in \emph{Ell}_{\fp}\\j(E)=0}}a(E)^n.$$
Now we may apply Theorem \ref{schoof's theorem} to obtain
$$\sum_{t=2}^{p-1}a(t,p)^n=\sum_{0<s<2\sqrt{p}}H(s^2-4p)s^n-\frac{1}{2}\sum_{\substack{[E]_{\fp}\in \emph{Ell}_{\fp}\\j(E)=1728}}a(E)^n-\frac{1}{2}\sum_{\substack{[E]_{\fp}\in \emph{Ell}_{\fp}\\j(E)=0}}a(E)^n.$$
Recall from (\ref{defn of H}) that if $d$ is the discriminant of an imaginary quadratic order $\mathcal{O}$,
$$H(d):=\sum_{\mathcal{O}\subseteq\mathcal{O}'\subseteq\mathcal{O}_{max}}h(\mathcal{O}'),$$ where the sum is over all orders between $\mathcal{O}$ and the maximal order.  Then taking $\ell$ as defined as in the statement of the Proposition, we have
$$H(s^2-4p)=\sum_{f|\ell} h\left(\frac{s^2-4p}{f^2}\right),$$
which gives  
\begin{equation}\label{h step}
\sum_{t=2}^{p-1}a(t,p)^n=\sum_{0<s<2\sqrt{p}}s^n\sum_{f|\ell} h\left(\frac{s^2-4p}{f^2}\right)-\frac{1}{2}\sum_{\substack{[E]_{\fp}\in \emph{Ell}_{\fp}\\j(E)=1728}}a(E)^n-\frac{1}{2}\sum_{\substack{[E]_{\fp}\in \emph{Ell}_{\fp}\\j(E)=0}}a(E)^n.
\end{equation}
To complete the proof, we apply Lemma \ref{h to h*} to the right side of \eqref{h step}, to replace $h$ by $h^*$. Then, collecting terms gives
\begin{align*}
\sum_{t=2}^{p-1}a(t,p)^n&=\sum_{0<s<2\sqrt{p}}s^n \sum_{f|\ell}h^*\left(\frac{s^2-4p}{f^2}\right)-\frac{1}{4}\sum_{\substack{[E]_{\fp}\in \emph{Ell}_{\fp}\\j(E)=1728}}a(E)^n-\frac{1}{6}\sum_{\substack{[E]_{\fp}\in \emph{Ell}_{\fp}\\j(E)=0}}a(E)^n\\
&=\sum_{0<s<2\sqrt{p}}s^n \sum_{f|\ell}h^*\left(\frac{s^2-4p}{f^2}\right)-2^{n-1}(a^n+b^n)\\
 &\quad-\frac{1}{3}[(c+d)^n+(2c-d)^n+(c-2d)^n],
\end{align*}
by Lemmas \ref{j inv 1728} and \ref{j inv 0}.  This is the desired result.
\end{proof}

%%%%%%%%%%%%%%%%%%%%%%%%
%% proof of theorem 2 %%
%%%%%%%%%%%%%%%%%%%%%%%%
Proposition \ref{step 1 of thm 2 proof} and Theorem \ref{hijikata} give us the tools necessary to complete the proof our second main theorem:
\begin{proof}[Proof of Theorem \ref{second main theorem}] By Theorem \ref{hijikata}, we have for $k\geq4$ even,
\allowdisplaybreaks{
\begin{align*}
\textnormal{Tr}_k(\Gamma,p)&=-1-\frac{1}{2}h^*(-4p)(-p)^{\frac{k}{2}-1}-\sum_{0<s<2\sqrt{p}}G_k(s,p)\sum_fh^*\left(\frac{s^2-4p}{f^2}\right)\\
                   &=-1-\frac{1}{2}h^*(-4p)(-p)^{\frac{k}{2}-1}-\sum_{0<s<2\sqrt{p}}\Biggl[(-p)^{\frac{k}{2}-1}\Biggr.\\
                   &\hspace*{0.25in}\Biggl.+\sum_{j=0}^{\frac{k}{2}-2}(-1)^j\binom{k-2-j}{j}p^js^{k-2j-2}\Biggr]\sum_fh^*\left(\frac{s^2-4p}{f^2}\right)\\
                   &=-1-\frac{1}{2}h^*(-4p)(-p)^{\frac{k}{2}-1}-(-p)^{\frac{k}{2}-1}\sum_{0<s<2\sqrt{p}}1\sum_fh^*\left(\frac{s^2-4p}{f^2}\right)\\
                   &\hspace*{0.2in}-\sum_{j=0}^{\frac{k}{2}-2}(-1)^j\binom{k-2-j}{j}p^j\sum_{0<s<2\sqrt{p}}s^{k-2j-2}\sum_fh^*\left(\frac{s^2-4p}{f^2}\right),
\end{align*}}after substituting in the definition of $G_k(s,p)$ and distributing.  Now, note that taking $k=2$ in Theorem \ref{hijikata} provides
$$0=p-\frac{1}{2}h^*(-4p)-\sum_{0<s<2\sqrt{p}}1\sum_fh^*\left(\frac{s^2-4p}{f^2}\right).$$
We apply this, together with Proposition \ref{step 1 of thm 2 proof} and obtain
\begin{align*}
\textnormal{Tr}_k(\Gamma,p)&=-1-\frac{1}{2}h^*(-4p)(-p)^{\frac{k}{2}-1}+(-p)^{\frac{k}{2}-1}\left(\frac{1}{2}h^*(-4p)-p\right)\\
                  &\hspace*{0.2in}-\sum_{j=0}^{\frac{k}{2}-2}(-1)^j\binom{k-2-j}{j}p^j\Biggl[\sum_{t=2}^{p-1}a(t,p)^{k-2j-2}+\frac{1}{2}\left[2^{k-2j-2}(a^{k-2j-2}\right. \Biggr.\\
                   &\left. \hspace*{0.2in} +b^{k-2j-2})\right] \Biggl.+\frac{1}{3}\left[(c+d)^{k-2j-2}+(2c-d)^{k-2j-2}+(c-2d)^{k-2j-2}\right]\Biggr]\\
                   &=-1+(-p)^{\frac{k}{2}-1}\cdot(-p)-\sum_{j=0}^{\frac{k}{2}-2}(-1)^j\binom{k-2-j}{j}p^j\sum_{t=2}^{p-1}a(t,p)^{k-2j-2}\\
                   &\hspace*{0.2in}-\frac{1}{2}\sum_{j=0}^{\frac{k}{2}-2}(-1)^j\binom{k-2-j}{j}p^j\left[(2a)^{k-2j-2}+(2b)^{k-2j-2}\right]\\
                   &\hspace*{0.2in}-\frac{1}{3}\sum_{j=0}^{\frac{k}{2}-2}(-1)^j\binom{k-2-j}{j}p^j\left[(c+d)^{k-2j-2}+(2c-d)^{k-2j-2}\right.\\&\hspace*{0.2in}\left.+(c-2d)^{k-2j-2}\right],
\end{align*}
after distributing once again.  Now, we notice the simple fact that
$$(-p)^{\frac{k}{2}-1}\cdot(-p)=-(-p)^{\frac{k}{2}-1}(p-2)-2\left(\frac{1}{2}(-p)^{\frac{k}{2}-1}\right)-3\left(\frac{1}{3}(-p)^{\frac{k}{2}-1}\right).$$
Splitting up the factors of $(-p)^{\frac{k}{2}-1}$ in this way gives that
\begin{align*}
\textnormal{Tr}_k(\Gamma,p) &=-1-\frac{1}{2}\left[G_k(2a,p)+G_k(2b,p)\right]\\
                            &\hspace*{0.2in}-\frac{1}{3}\left[G_k(c+d,p)+G_k(2c-d,p)+G_k(c-2d,p)\right]\\
               &\hspace{0.2in}-(p-2)(-p)^{\frac{k}{2}-1}-\sum_{j=0}^{\frac{k}{2}-2}(-1)^j\binom{k-2-j}{j}p^j\sum_{t=2}^{p-1}a(t,p)^{k-2j-2}\\
               &=-1-\lambda(k,p)-\sum_{t=2}^{p-1}(-p)^{\frac{k}{2}-1}-\sum_{t=2}^{p-1}\sum_{j=0}^{\frac{k}{2}-2}(-1)^j\binom{k-2-j}{j}p^ja(t,p)^{k-2j-2}\\
               &=-1-\lambda(k,p)-\sum_{t=2}^{p-1}G_k(a(t,p),p),
\end{align*}
according to the definitions of $G_k$ and $\lambda(k,p)$ given in the statement of the theorem.  This completes the proof of Theorem \ref{second main theorem}.
\end{proof} 

\begin{remark}\label{remk}
According to Theorem \ref{first main theorem}, we may rewrite $a(t,p)$ in terms of the hypergeometric function $\hg{\xi}{\xi^5}{\varepsilon}{t}$.  Thus, Theorem \ref{second main theorem} can be reformulated to give $\textnormal{Tr}_k(\Gamma,p)$ in terms of $\lambda(k,p)$ and $G_k\left(\psi^{-1}(t)p\,\hg{\xi}{\xi^5}{\varepsilon}{t},p\right)$, where $\psi(t)=-\phi(2)\xi^{-3}(1-t).$
\end{remark}

%%%%%%%%%%%%%%%%%%%%%%%%%
% Proof of Recursion %%%%
%%%%%%%%%%%%%%%%%%%%%%%%%

\section{Proof of Theorem \ref{recursion}}\label{recursion proof}

Theorem \ref{recursion} is proved by combining Theorems \ref{first main theorem} and \ref{second main theorem} with an inverse pair given in \cite{Ri68}.  First, recall that 

$$G_k(s,p)=\sum_{j=0}^{\frac{k}{2}-1} (-1)^j \binom{k-2-j}{j} p^j s^{k-2j-2}.$$
Letting $m=\frac{k}{2}-1$ and $H_m(x) :=\sum_{i=0}^{m}\binom{m+i}{m-i}x^i$, we have that
\begin{equation}\label{G and H}
G_k(s,p)=(-p)^mH_m\left(\frac{-s^2}{p}\right).
\end{equation}
Now, we make use of the inverse pair \cite[p. 67]{Ri68} given by
\begin{equation}\label{inverse pair}
\rho_n(x)=\sum_{k=0}^{n}\binom{n+k}{n-k}x^k, \qquad x^n=\sum_{k=0}^n(-1)^{k+n}\left[\binom{2n}{n-k}-\binom{2n}{n-k-1}\right]\rho_k(x).
\end{equation}
Applied to the definition of $H_m$, this gives
\begin{equation*}
x^m=\sum_{i=0}^{m}(-1)^{i+m}\left[\binom{2m}{m-i}-\binom{2m}{m-i-1}\right]H_i(x).
\end{equation*}
By taking $x=\frac{-s^2}{p}$, together with \eqref{G and H}, we have

\begin{eqnarray}\label{s2m eqn}
s^{2m}&=&\sum_{i=0}^{m}p^{m-i}\left[\binom{2m}{m-i}-\binom{2m}{m-i-1}\right]G_{2i+2}(s,p)\notag\\
&=&\sum_{i=0}^{m}b_iG_{2i+2}(s,p), 
\end{eqnarray}
where $b_i$ is as defined in the statement of Theorem \ref{recursion}.  

\begin{proof}[Proof of Theorem \ref{recursion}]  By \eqref{s2m eqn}, we have 
\begin{eqnarray}\label{s2m eqn2}
s^{2m}&=&\sum_{i=0}^{m}b_iG_{2i+2}(s,p)\notag \\
&=&G_{2m+2}(s,p)+\sum_{i=0}^{m-1}b_iG_{2i+2}(s,p),
\end{eqnarray}
since $b_m=1$.  Now, for $m\geq1$, Theorem \ref{second main theorem} implies
\begin{align*}
\textnormal{Tr}_{2(m+1)}(\Gamma,p)&= -1-\lambda(2m+2,p)-\sum_{t=2}^{p-1}G_{2m+2}(a(t,p),p)\\
&=-1-\lambda(2m+2,p)-\sum_{t=2}^{p-1}\Biggl( a(t,p)^{2m} \Biggr. \\
&\left. \hspace*{.2in}-\sum_{i=0}^{m-1}b_iG_{2i+2}(a(t,p),p)\right) \qquad (\textnormal{by \eqref{s2m eqn2}})\\
%&=&-1-\lambda(2m+2,p)-\sum_{t=2}^{p-1}a(t,p)^{2m}+\sum_{i=0}^{m-1}b_i\sum_{t=2}^{p-1}G_{2i+2}(a(t,p))\\
&=-1-\lambda(2m+2,p)-\sum_{t=2}^{p-1}a(t,p)^{2m}+b_0\sum_{t=2}^{p-1}G_2(a(t,p),p)\\
&	\hspace*{0.2in}+\sum_{i=1}^{m-1}b_i\sum_{t=2}^{p-1}G_{2i+2}(a(t,p),p)\\
&=-1-\lambda(2m+2,p)-\sum_{t=2}^{p-1}a(t,p)^{2m}+b_0(p-2)\\
&\hspace*{0.2in}-\sum_{i=1}^{m-1}b_i(\textnormal{Tr}_{2i+2}(\Gamma,p)+1+\lambda(2i+2,p)),
\end{align*} 
by Theorem \ref{second main theorem} and since $G_2=1$.  The proof is completed by rearranging and noting that Theorem \ref{first main theorem} implies

\begin{align*}
a(t,p)^{2m}&=p^{2m}\phi^{2m}(2)\xi^{6m}(1-t)\hg{\xi}{\xi^5}{\varepsilon}{t}^{2m}\\
&=p^{2m}\phi^m(1-t)\hg{\xi}{\xi^5}{\varepsilon}{t}^{2m},
\end{align*}
since $\phi^2=\varepsilon$ and $\xi^6=\phi$.
\end{proof}

%%%%%%%%%%%%%%%%%%%%%%%%%%
%% tau corollaries %%%%%%%
%%%%%%%%%%%%%%%%%%%%%%%%%%

\section{$\tau(p)$ Corollaries}\label{tau corollaries}

Specializing to various values of $k$ in Theorem \ref{second main theorem}, we arrive at more explicit formulas.  In particular, by taking $k=12$ we obtain a formula for Ramanujan's $\tau$-function.  Recall that we define $\tau(n)$ by
$$(2\pi)^{-12}\Delta(z)=q\prod_{n=1}^{\infty}(1-q^n)^{24}=\sum_{n=1}^{\infty}\tau(n)q^n.$$
Also, recall that $\Delta(z)$ generates the one dimensional space $S_{12}$, and thus $\textnormal{Tr}_{12}(\Gamma,p)=\tau(p)$ for primes $p$.  We conclude with results that stem from this specialization to $k=12$.  

Throughout, we let $a,b\in\mathbb{Z}$ be such that  $p=a^2+b^2$ and $a+bi \equiv 1 \,(2+2i)$ in $\mathbb{Z}[i]$.  Further, let $c,d\in\mathbb{Z}$ satisfy $p=c^2-cd+d^2$ and $c+d\omega \equiv 2 \,(3)$ in $\mathbb{Z}[\omega]$, where $\omega=e^{2\pi i /3}$.  The first corollary follows from a straightforward application of Theorem \ref{second main theorem}:

\begin{corollary}\label{tau corollary 1}
Let  $a,b,c,$ and $d$ be defined as above, and set $x=a^2b^2$ 
and $y=cd$.  If $p$ is a prime, $p \equiv 1 (12)$, then
$$\tau(p)=-1-8p^5+80p^3x-256px^2+27y^2p^3-27y^3p^2 
-\sum_{t=2}^{p-1} G_{12}(a(t,p),p),$$
where
$$G_{12}(s,p)=s^{10}-9ps^8+28p^2s^6-35p^3s^4+15p^4s^2-p^5.$$
\end{corollary}

As noted previously, \ref{tau corollary 1} can be reformulated in terms of the hypergeometric function $\hg{\xi}{\xi^5}{\varepsilon}{t}$.  In fact, we can inductively arrive at a formula for $\tau(p)$ in terms only of $10^{th}$ powers of this hypergeometric function.

\begin{corollary}\label{tau corollary 2}
Let $p\equiv 1 \pmod{12}$ be prime and let $a$, $b$, $c$, and $d$ be defined as above.  Let $\xi$ be an element of order $12$ in $\widehat{\mathbb{F}_p^{\times}}$.  Then
\begin{align*}
\tau(p)&=42p^6-90p^4-75p^3-35p^2-9p-1-2^9(a^{10}+b^{10})\\
&\quad -\frac{1}{3}\left( (c+d)^{10}+(2c-d)^{10}+(c-2d)^{10} \right) - \sum_{t=2}^{p-1} p^{10}\phi(1-t)\hg{\xi}{\xi^5}{\varepsilon}{t}^{10}.
\end{align*}
\end{corollary}

\begin{proof}
Recall from the statement of Theorem \ref{second main theorem} that we have
$$\lambda(k,p)=\frac{1}{2}[G_k(2a,p)+G_k(2b,p)]+\frac{1}{3}[G_k(c+d,p)+G_k(2c-d,p)+G_k(c-2d,p)],$$
where 
$$G_k(s,p)=\sum_{j=0}^{\frac{k}{2}-1}(-1)^j\binom{k-2-j}{j}p^js^{k-2j-2}.$$
Then, one can check by hand or with Maple that we have the following, recalling the relations $p=a^2+b^2=c^2-cd+d^2$:
\begin{align*}\label{lambdas}
\lambda(4,p)&=2p\\
\lambda(6,p)&=-4p^2+2^3(a^4+b^4)\\
\lambda(8,p)&=-8p^3+2^5(a^6+b^6)-40p(a^4+b^4)+\frac{1}{3}((c+d)^6+(2c-d)^6+(c-2d)^6)\\
\lambda(10,p)&=52p^4+2^7(a^8+b^8)-224p(a^6+b^6)+120p^2(a^4+b^4)\\
&\quad+\frac{1}{3}((c+d)^8+(2c-d)^8+(c-2d)^8)-\frac{7}{3}p((c+d)^6+(2c-d)^6+(c-2d)^6)\\
\lambda(12,p)&=-152p^5+2^{9}(a^{10}+b^{10})-1152p(a^8+b^8)+896p^2(a^6+b^6)-280p^3(a^4+b^4)\\
&\quad+\frac{1}{3}((c+d)^{10}+(2c-d)^{10}+(c-2d)^{10})\\&\quad-3p((c+d)^8+(2c-d)^8+(c-2d)^8)\\
&\quad+\frac{28}{3}p^2((c+d)^6+(2c-d)^6+(c-2d)^6).
\end{align*}
By using each successive formula for $G_k(a(t,p),p)$, together with Theorem \ref{second main theorem} and the formulas for $\lambda(k,p)$ given above, we can compute $\sum_{t=2}^{p-1}a(t,p)^{k-2}$, for $k=4,\dots,12$.  We exhibit the computations in the cases $k=4$ and $k=6$, to give the idea of the technique.  First, notice that $G_4(s,p)=s^2-p$ and recall $\textnormal{Tr}_4(\Gamma,p)=0$, as there are no cusp forms of weight $4$ for $\Gamma$.  Thus, Theorem \ref{second main theorem} implies
\begin{align*}
0=\textnormal{Tr}_4(\Gamma,p)&=-1-\lambda(4,p)-\sum_{t=2}^{p-1}G_4(a(t,p),p)\\
&=-1-2p-\sum_{t=2}^{p-1}(a(t,p)^2-p)\\
&=-1-2p+p(p-2)-\sum_{t=2}^{p-1}a(t,p)^2.
\end{align*}
Thus, after simplifying, we see that
\begin{equation}\label{power 2}
0=p^2-4p-1-\sum_{t=2}^{p-1}a(t,p)^2.
\end{equation}
Now, we utilize this computation to derive a formula for the sum of $4^{th}$ powers of $a(t,p)$.  For $k=6$, we have $G_6(s,p)=s^4-3ps^2+p^2$ and once again  $\textnormal{Tr}_6(\Gamma,p)=0$.  Then, by Theorem \ref{second main theorem} and the formula for $\lambda(6,p)$ given at the start of the proof, we see that
\begin{align*}
0=\textnormal{Tr}_6(\Gamma,p)&=-1-\lambda(6,p)-\sum_{t=2}^{p-1}G_6(a(t,p),p)\\
&=-1+4p^2-2^3(a^4+b^4)-\sum_{t=2}^{p-1}(a(t,p)^4-3pa(t,p)^2+p^2).
\end{align*}
We distribute the summation across the polynomial $G_6(a(t,p),p)$ and then make a substitution for $\sum_{t=2}^{p-1}a(t,p)^2$, according to (\ref{power 2}).  This gives
\begin{align*}
0=\textnormal{Tr}_6(\Gamma,p)&=4p^2-1-2^3(a^4+b^4)-\sum_{t=2}^{p-1}a(t,p)^4+3p\sum_{t=2}^{p-1}a(t,p)^2-p^2(p-2)\\
&=4p^2-1-2^3(a^4+b^4)-\sum_{t=2}^{p-1}a(t,p)^4\\
&\hspace*{0.75in}+3p(-2p-1+p(p-2))-p^2(p-2).
\end{align*}
After simplifying, we arrive at
\begin{equation}\label{power 4}
0=2p^3-6p^2-3p-2^3(a^4+b^4)-1-\sum_{t=2}^{p-1}a(t,p)^4.
\end{equation}
We continue this process, using successive formulas for $G_k(s,p)$ and $\lambda(k,p)$ and back-substituting previous results such as (\ref{power 2}) and (\ref{power 4}).  We omit the tedious details of the next couple of cases, which result in the following:
\begin{multline}\label{power 6}
\textnormal{Tr}_8(\Gamma,p)=0=5p^4-9p^2-5p-1-2^5(a^6+b^6)\\-\frac{1}{3}((c+d)^6+(2c-d)^6+(c-2d)^6)-\sum_{t=2}^{p-1}a(t,p)^6
\end{multline}
\begin{multline}\label{power 8}
\textnormal{Tr}_{10}(\Gamma,p)=0=14p^5-28p^3-20p^2-7p-1-2^7(a^8+b^8)\\-\frac{1}{3}((c+d)^8+(2c-d)^8+(c-2d)^8)-\sum_{t=2}^{p-1}a(t,p)^8
\end{multline}

Now, we use (\ref{power 2}),$\dots$,(\ref{power 8}) together with the formula for $\lambda(12,p)$ from the beginning of the proof to compute a formula for $\tau(p)$.  Since $G_{12}(s,p)=s^{10}-9ps^8+28p^2s^6-35p^3s^4+15p^4s^2-p^5$ and $\textnormal{Tr}_{12}(\Gamma,p)=\tau(p)$, Theorem \ref{second main theorem} gives
\begin{align*}
\tau(p)&=-1-\lambda(12,p)-\sum_{t=2}^{p-1}G_{12}(a(t,p),p)\\
&=-1-\lambda(12,p)-\sum_{t=2}^{p-1}a(t,p)^{10}+9p\sum_{t=2}^{p-1}a(t,p)^8-28p^2\sum_{t=2}^{p-1}a(t,p)^6\\
&\hspace*{0.2in}+35p^3\sum_{t=2}^{p-1}a(t,p)^4-15p^4\sum_{t=2}^{p-1}a(t,p)^2+p^5(p-2)
\end{align*}
%&=-1-\lambda(12,p)-\sum_{t=2}^{p-1}a(t,p)^{10}\\
%&\hspace*{0.2in}+9p\Bigl(14p^5-28p^3-20p^2-7p-1-2^7(a^8+b^8)\Bigr.\\
%&\hspace*{0.2in}\Bigl.-\frac{1}{3}((c+d)^8+(2c-d)^8+(c-2d)^8)\Bigr)\\
%&\hspace*{0.2in}-28p^2\Bigl(5p^4-9p^2-5p-1-2^5(a^6+b^6)-\frac{1}{3}((c+d)^6+(2c-d)^6+(c-2d)^6)\Bigr)\\
%&\hspace*{0.2in}+35p^3(2p^3-6p^2-3p-2^3(a^4+b^4)-1)-15p^4(p^2-4p-1)+p^5(p-2),
%\end{align*}
%by (\ref{power 2}), $\dots$,(\ref{power 8}). The substitution
%\begin{align*}
%\lambda(12,p)&=-152p^5+2^{9}(a^{10}+b^{10})-1152p(a^8+b^8)+896p^2(a^6+b^6)-280p^3(a^4+b^4)\\
%&\quad+\frac{1}{3}((c+d)^{10}+(2c-d)^{10}+(c-2d)^{10})\\&\quad-3p((c+d)^8+(2c-d)^8+(c-2d)^8)\\
%&\quad+\frac{28}{3}p^2((c+d)^6+(2c-d)^6+(c-2d)^6).

Substitutions via (\ref{power 2}),$\dots$,(\ref{power 8}), together with the formula for $\lambda(12,p)$, give rise to many cancellations, resulting in 
\begin{multline}\label{power 10}
\tau(p)=42p^6-90p^4-75p^3-35p^2-9p-1-2^9(a^{10}+b^{10})\\-\frac{1}{3}((c+d)^{10}+(2c-d)^{10}+(c-2d)^{10})-\sum_{t=2}^{p-1}a(t,p)^{10}.
%\end{align*}
\end{multline}
%gives rise to many cancellations, and after simplifying, we find that
%\begin{multline}\label{power 10}
%\tau(p)=42p^6-90p^4-75p^3-35p^2-9p-1-2^9(a^{10}+b^{10})\\-\frac{1}{3}((c+d)^{10}+(2c-d)^{10}+(c-2d)^{10})-\sum_{t=2}^{p-1}a(t,p)^{10}.
%\end{multline}

 Finally, to complete the proof, we recall that Theorem \ref{first main theorem} implies 
$$p\, \hg{\xi}{\xi^5}{\varepsilon}{t}=-\phi(2)\xi^{-3}(1-t)a(t,p)$$
for $t\in\fp \backslash \{0,1\}$.  Thus, 
$$a(t,p)^{10}=\left(-p\phi(2)\xi^3(1-t)\hg{\xi}{\xi^5}{\varepsilon}{t}\right)^{10}=p^{10}\phi(1-t)\hg{\xi}{\xi^5}{\varepsilon}{t}^{10},$$
since $\phi^{10}=\varepsilon$ and $\xi^{30}=\xi^6=\phi$.  This, together with (\ref{power 10}), confirms the corollary.
\end{proof}

The technique used in the proof of Corollary \ref{tau corollary 2} can be extended one step further to arrive at yet another formula for $\tau(p)$:
\begin{corollary}\label{tau corollary 3}
Let $p\equiv 1 \pmod{12}$ be prime and let $a$, $b$, $c$, and $d$ be defined as above.  Let $\xi$ be an element of order $12$ in $\widehat{\mathbb{F}_p^{\times}}$.  Then
\begin{align*}
\tau(p)&=12p^6-27p^4-25p^3-14p^2-\frac{54}{11}p-1-\frac{1}{11p}-\frac{2^{11}}{11p}(a^{12}+b^{12})\\
&\hspace*{0.2in}-\frac{1}{33p}((c+d)^{12}+(2c-d)^{12}+(c-2d)^{12})-\frac{1}{11}\sum_{t=2}^{p-1}p^{11}\,\hg{\xi}{\xi^5}{\varepsilon}{t}^{12}.
\end{align*}
\end{corollary}

\begin{proof}
The proof follows in the same way as for Corollary \ref{tau corollary 2}, but by taking $k=14$, so we give a sketch here.  First, one calculates that 
%We continue the process used in the proof of Corollary \ref{tau corollary 2} to consider the case $k=14$.  Using the notation from the statement of Theorem \ref{second main theorem}, one calculates that
\begin{equation*}
G_{14}(s,p)=s^{12}-11ps^{10}+45p^2s^8-84p^3s^6+70p^4s^4-21p^5s^2+p^6.
\end{equation*}
Now, since there are no cusp forms of level $14$ for $\Gamma$, we have $\textnormal{Tr}_{14}(\Gamma,p)=0$.  Thus, Theorem \ref{second main theorem} implies
\begin{equation*}
0=-1-\lambda(14,p)-\sum_{t=2}^{p-1}G_{14}(a(t,p),p).
\end{equation*}
Now, by applying (\ref{power 2}),$\dots$,(\ref{power 10}), together with a formula for $\lambda(14,p)$, given by
%&=-1-\lambda(14,p)-\sum_{t=2}^{p-1}a(t,p)^{12}+11p\sum_{t=2}^{p-1}a(t,p)^{10}-45p^2\sum_{t=2}^{p-1}a(t,p)^8\\
%&\hspace*{0.2in}+84p^3\sum_{t=2}^{p-1}a(t,p)^6-70p^4\sum_{t=2}^{p-1}a(t,p)^4+21p^5\sum_{t=2}^{p-1}a(t,p)^2-p^6(p-2).
%\end{align*}
%Next, we make a substitution for each $\displaystyle{\sum_{t=2}^{p-1}a(t,p)^n}$ for $n=2,\dots,10$, according to (\ref{power 2}),$\dots$,(\ref{power 10}).  Then we have
%\begin{align*}
%0&=-1-\lambda(14,p)-\sum_{t=2}^{p-1}a(t,p)^{12}+11p\Bigl(-\tau(p)+42p^6-90p^4-75p^3-35p^2-9p-1\Bigr.\\
%&\hspace*{0.2in}\Bigl.-2^9(a^{10}+b^{10})-\frac{1}{3}((c+d)^{10}+(2c-d)^{10}+(c-2d)^{10})\Bigr)\\
%&\hspace*{0.2in}-45p^2\Bigl(14p^5-28p^3-20p^2-7p-1-2^7(a^8+b^8)\\
%&\hspace*{0.2in}-\frac{1}{3}((c+d)^8+(2c-d)^8+(c-2d)^8)\Bigr)\\
%&\hspace*{0.2in}+84p^3\Bigl(5p^4-9p^2-5p-1-2^5(a^6+b^6)-\frac{1}{3}((c+d)^6+(2c-d)^6+(c-2d)^6)\Bigr)\\
%&\hspace*{0.2in}-70p^4(2p^3-6p^2-3p-1)+21p^5(p^2-4p-1)-p^6(p-2).
%\end{align*}
%One may calculate, by hand or with the aide of Maple, that we have
\begin{align*}
  \lambda(14,p)&=338p^6+2^{11}(a^{12}+b^{12})-11\cdot2^9p(a^{10}+b^{10})+45\cdot2^7p^2(a^8+b^8)\\
&\hspace*{0.2in}-84\cdot2^5p^3(a^6+b^6)+70\cdot2^3p^4(a^4+b^4)\\
&\hspace*{0.2in}+\frac{1}{3}((c+d)^{12}+(2c-d)^{12}+(c-2d)^{12})\\
&\hspace*{0.2in}-\frac{11}{3}p((c+d)^{10}+(2c-d)^{10}+(c-2d)^{10})\\
&\hspace*{0.2in}+15p^2((c+d)^8+(2c-d)^8+(c-2d)^8)\\
&\hspace*{0.2in}-28p^3((c+d)^6+(2c-d)^6+(c-2d)^6),
\end{align*}
one makes cancellations and finds that
\begin{multline}\label{power 12}
0=132p^7-297p^5-275p^4-154p^3-54p^2-1-11p-11p\tau(p)-2^{11}(a^{12}+b^{12})\\-\frac{1}{3}((c+d)^{12}+(2c-d)^{12}+(c-2d)^{12})-\sum_{t=2}^{p-1}a(t,p)^{12}.
\end{multline}
We now recall that Theorem \ref{first main theorem} implies 
$$p\, \hg{\xi}{\xi^5}{\varepsilon}{t}=-\phi(2)\xi^{-3}(1-t)a(t,p)$$
for $t\in\fp \backslash \{0,1\}$.  Therefore,
$$a(t,p)^{12}=\left(-p\phi(2)\xi^3(1-t)\hg{\xi}{\xi^5}{\varepsilon}{t}\right)^{12}=p^{12}\, \hg{\xi}{\xi^5}{\varepsilon}{t}^{12},$$
since $\phi$ has order $2$ and $\xi$ has order $12$.  Making this substitution into (\ref{power 12}) and then solving for $\tau(p)$ gives the desired result.

\end{proof}

\section*{Acknowledgments}
The author thanks her advisor M. Papanikolas for his advice and support during the preparation of this paper.  The author also thanks the Department of Mathematics at Texas A$\&$M University, where the majority of this research was conducted.

\end{document}